# MOUVEMENT
# &
# ORIGINE DU CALCUL INFINITÉSIMAL

## Philosophie et continuité

Salomon OFMAN
Institut mathématique de Jussieu-Paris Rive Gauche
Histoire des Sciences mathématiques
4 Place Jussieu
75005 Paris
salomon.ofman@imj-prg.fr

'*Motion consists merely in the fact that bodies are sometimes in one place and sometimes in another, and that they are at intermediate places at intermediate times. Only those who have waded through the quagmire of philosophic speculation on this subject can realise what a liberation from antique prejudices is involved in this simple and straightforward commonplace*' (B. Russell, *Mathematics and Metaphysicians*, p. 66).

*English Abstract.* – The text has his origin in some lectures given at the University of Bologna, inside an interdisciplinary program of mathematics, history of science, physics and philosophy. Since they are at the junction of these fields, movement and infinitesimal calculus are good instances to understand some fundamental questions in history of sciences, the diversity of the conceptions stretching across it and more generally the philosophy background of theses researches. It is the field where the 'struggle of giants' takes place.
             In this article we consider the general foundations needed for this study, in particular its philosophical background. While the principal question concerns 'the one and the multiple' it is treated through the study of the 'discreet and the continuous' considered under both his philosophical and mathematical aspect.
             In a second article ([OFM2]) we will consider the relation between the theorization of the movement and the infinitesimal calculus, and their historical background where they were able to grow.

Résumé. – Ce texte a son origine dans des exposés donnés dans le cadre d'un programme de philomathique, c'est-à-dire trans-disciplinaire, mathématiques, histoire des sciences, physique et philosophie. Le mouvement, mais aussi le calcul infinitésimal, qui sont à l'intersection de ces disciplines, paraissent un bon exemple pour comprendre certaines questions fondamentales en histoire des sciences, les conceptions qui la traversent, et plus généralement le cadre philosophique dans lequel s'inscrivent ces recherches. C'est le champ où se déroule une 'bataille de géants'.

Dans ce premier article, nous posons les fondements généraux nécessaires à cette étude, en particulier son contexte philosophique. Si la question principale est celle de l'un et du multiple, son traitement passe par celle du discret et du continu, sous le double aspect mathématique et philosophique.

Cet éclairage conduit, dans un second texte ([OFM2]), à interpréter de manière quelque peu différente des points de vue classiques, le rapport entre théorisation du mouvement et calcul infinitésimal, ainsi que la situation historique dans laquelle ils ont pu se développer.

La référence aux textes de Galilée renvoie au tome et à la page de l'Édition nationale italienne.

**I. Introduction.**

'οὐ λόγους, ἀλλ' ὃ ὑμεῖς τιμᾶτε [ὦ ἄνδρες 'Αθηναῖοι], ἔργα' ('non pas des mots, mais ce que vous [hommes d'Athènes] estimez le plus, des actes') (Platon, *Apologie de Socrate*, 32A4-5).

Au travers de la théorisation du mouvement et du calcul infinitésimal, il s'agit de comprendre l'articulation entre le flux continu de la connaissance et la rupture de la découverte en sciences, mais aussi la diversité des controverses à leur sujet. Ces questions seront considérées sous un quadruple aspect.

Tout d'abord, nous voudrions montrer que la multiplicité des interprétations portant sur le mouvement et les infinitésimaux, trouve son origine dans des grands courants de pensée qui, loin de se borner à l'histoire des idées ou même aux disciplines dites classiques, traversent le champ scientifique tout entier, que ce soit en physique en biologie ou en paléontologie ou, plus étonnamment encore, en mathématiques. Ces multiples courants s'ordonnent entre deux positions extrêmes qui se ramènent à une alternative concernant la conception de l'histoire et de la philosophie des sciences. Celle-ci, à son tour, apparaît comme rien moins qu'un choix entre deux conceptions du monde.

Cette réflexion sur la pensée scientifique et son histoire, se fonde donc sur une approche philosophique, en entendant par « philosophie » une démarche unificatrice suffisamment large, en quelque sorte une universalisation qui nécessairement dépasse les frontières usuelles entre disciplines. En tant que recherche sur les origines, il faut solliciter les écrits les plus anciens, les ouvrages de la Grèce ancienne, principalement ceux de Platon et d'Aristote, que nous confronterons à des textes plus récents, principalement ceux de Poincaré.

Nous considérons alors comment notre sujet et les différentes approches des historiens et des philosophes des sciences s'insèrent dans cette problématique générale.

Dans [OFM2], nous analysons les notions aristotéliciennes sur le mouvement, en tant qu'elles prédominent jusqu'à Galilée. L'étude des textes de Platon, d'Aristote et surtout d'Euclide conduit à une critique de certaines interprétations usuelles. Nous examinons ensuite l'exposition par Galilée de la science du mouvement. C'est dans le *Dialogue sur les deux grands systèmes du monde* que l'on peut mesurer les difficultés auxquelles son auteur a été confronté dans la construction de cette science. La dernière partie de [OFM2] porte sur les conséquences de cette nouvelle science et sa relation aux infinitésimaux. En conclusion, nous nous proposons inversement de comprendre l'impact de la théorie des infinitésimaux sur la conception du mouvement, ainsi que les causes qui ont pu empêcher les mathématiciens/physiciens antérieurs d'élaborer une théorie du mouvement, au sens où nous l'entendons.

## II. Présentation

### 1. Sur le calcul infinitésimal et le mouvement.

Le calcul infinitésimal a amené des changements fondamentaux non seulement dans toutes les branches des mathématiques, mais aussi dans la plupart des domaines de la science. C'est 'l'un des plus formidables accomplissements de l'esprit humain' (Richard Courant, cf. [BOY], avant-propos), 'un si puissant instrument, qu'il a changé la face des mathématiques pures et appliquées' ([COU], p. 201), 'l'instrument le plus effectif pour l'étude que les mathématiques aient jamais produit' (Carl Boyer, *Encyclopædia Britannica*, entrée : *Calculus*), 'la plus puissante méthode mathématique jamais inventée' ([SEI]).

Nous nous intéresserons ici au calcul infinitésimal en relation à l'étude du mouvement et de la vitesse. Mouvement et vitesse n'étant pas seulement des objets de la physique, mais aussi bien des sciences naturelles ou historiques que de la philosophie. Nous sommes ainsi conduit aux débuts de la physique moderne, qui naît avec la théorie mathématique de Galilée, présenté comme le 'père de la science moderne' ([CLA], préface). De même que l'introduction des infinitésimaux a bouleversé la totalité des mathématiques, la physique galiléenne a changé la physique et au-delà (cf. par exemple [KOY], p. 166). Cette étude se conclut dans [OFM2] sur le calcul infinitésimal, en le reliant historiquement à la question du mouvement/vitesse.

L'une des difficultés réside dans les mots eux-mêmes, ceux de 'mouvement', 'vitesse' mais aussi les termes convenables pour désigner ce qui est défini par les mathématiques. C'est en effet avec la science galiléenne que la vitesse acquiert un (nouveau) sens, d'où la nécessité d'éviter les pièges du langage, souvent perpétués par les historiens.

La question du mouvement a été l'un des problèmes les plus complexes et les plus étudiés. Pour les Grecs anciens, cette question est centrale et pose d'emblée un paradoxe à la pensée. Un mouvement est mouvement de quelque chose, mais ce quelque chose, pour être, doit, au moins d'une certaine manière, être au repos. Chez Platon, le mouvement est un point crucial pour unifier deux mondes, l'un éternel et immobile, l'autre temporel et toujours différent. Avec le repos, ce sont deux des cinq genres fondamentaux de l'être (*Sophiste*, 254d4-5). Pour Aristote, qui propose une solution aux paradoxes de Zénon, c'est par lui que la nature, conçue comme ce qui est sujet au mouvement, est définie. C'est donc l'objet essentiel de toute science physique (*Physique*, II, 192b14) et comme tel, il se trouve encore au cœur même de la philosophie.

### 2. A propos d'une question de terminologies.

De manière très générale, les approches de la délicate question des origines sont extrêmement diverses. Elles se situent entre deux positions qui dépassent de loin la seule histoire des sciences. Ainsi, les exemples actuels les plus marquants se trouvent en géologie. Ces positions, à leur tour, sont l'expression de deux camps philosophiques antagoniques, chacun alignant ses champions. Il s'agit ici d'Aristote pour l'un, de Galilée pour l'autre.

Notre manière de penser le mouvement et la vitesse est tributaire de la dynamique qui marque la naissance de la physique moderne, au sortir de la pensée scolastique, toute imprégnée de philosophie aristotélicienne (ou du moins des interprétations qu'elle s'en donne).

Pour nous, la question fondamentale de la théorisation aristotélicienne du mouvement, essentiellement dans la *Physique*, est son rapport au concept moderne de mouvement. Cela pose aussitôt celle de la vitesse et donc du découpage instantané du temps. Il s'ensuit que le point de vue moderne de vitesse ne pouvait faire conceptuellement sens, tant qu'elle n'était pas définie physiquement en tant qu'objet de mesure.

Beaucoup d'analyses s'avèrent alors non fondées, puisque d'une manière ou d'une autre, elles procèdent d'une certaine confusion entre différents sens (modernes) de vitesse. Pour éviter cette confusion, nous lui substituerons le terme de 'célérité', pour tout ce qui précède la physique galiléenne (il en est de même pour M. Caveing, qui préfère toutefois le terme de 'vélocité', cf. [CAV], p. 168).

La physique aristotélicienne apparaît alors comme l'aboutissement des connaissances physico-mathématiques de l'antiquité grecque, et non une sorte d'élaboration purement empirique dont la conséquence essentielle aurait été de freiner les progrès scientifiques.

Les travaux réalisés au Moyen-Âge en mathématiques concernant des calculs, qu'en termes modernes on qualifierait d'intégraux ou de différentiels, conduisent certains historiens à questionner le moment où naît cette science moderne (pourquoi Galilée), ou à s'interroger sur les raisons du 'retard' qu'elle aurait eu à émerger. Une même réponse paraît s'imposer : la nécessité d'une 'vitesse' mesurable. Mais en définitive, la difficulté (l'impossibilité) à penser '**la** vitesse', est celle de penser son unité. D'où le rôle des travaux galiléens sur la chute des corps, dont le résultat essentiel et paradoxal, est de mesurer quelque chose qui n'existait pas. Ce qui est mesuré peut alors être nommé justement 'vitesse'. Travaillée par les mathématiciens, elle ouvrait la voie à la formalisation des (calculs) infinitésimaux.

### 3. La position du problème.

Depuis l'antiquité, les philosophes, physiciens ou les métaphysiciens, se sont confrontés au mouvement. Mais s'accordaient-ils sur le sens de ce terme ou celui-ci cachait-il une équivoque ? Était-ce bien la même chose qui était définie ?

Le calcul infinitésimal, s'il ne pose pas la même question d'unicité de définition, n'en présente pas moins un problème d'origine. Celle-ci a donné lieu à l'une des plus longues polémiques, pour déterminer celui qui, le premier, a formalisé sa définition. Un historien a pu même intituler un livre sur la querelle entre Newton et Leibniz : 'des philosophes en guerre ('*Philosophers at war*') ([HAL]). Cette querelle n'est pas sans intérêt, car elle éclaire certains aspects concernant la communauté des scientifiques, mathématiciens ou physiciens, leur psychologie et celle des nations (par exemple, [WES], chap. 14). Et dans une certaine mesure, on la retrouve à propos de la relativité. Mais elle est extrinsèque à la présente recherche.

Elle recouvre toutefois une question importante et toujours très controversée. Les concepts mathématiques sont-ils endogènes, dérivant essentiellement des problèmes que les mathématiciens se posent, ou leurs origines sont-elles extérieures aux mathématiques ? D'autre part, est-il légitime de lier le mouvement et le calcul différentiel ?

Notre approche du mouvement met à jour une dualité à l'œuvre en histoire des sciences, mais aussi dans les théories scientifiques elles-mêmes, qui se clôt sur l'un des conflits majeurs de la pensée, un combat de géants dont parle Platon (cf. §III.4).

Pour le philosophe athénien, il est essentiel de dépasser cette alternative traditionnelle (éternelle ?). C'est la démarche que nous essaierons de suivre ici et on ne s'étonnera pas que nos conclusions puissent sembler paradoxales, aboutissant à une double 'réhabilitation', celle d'Aristote (surtout) et celle de Galilée (qui en a moins besoin), et par delà, à une certaine unification des problématiques chez ces deux auteurs.

Lorsqu'une théorie bouleverse la manière même de concevoir certaines questions fondamentales, deux questions se posent naturellement : 'comment y est-on parvenu ?', mais aussi sa symétrique, 'pourquoi n'y pas être parvenu auparavant ?'. Les réponses peuvent, là encore, être ordonnées comme des variations entre deux conceptions extrêmes. La plus commune considère quelque blocage 'psychologique' et/ou 'sociologique' ; l'autre, plus 'charitable', se fonde sur l'absence de certains outils 'techniques', laissant pendante la caractérisation de ce qui était absent.

Dans l'étude de ces problèmes, nous adopterons une démarche de généralité croissante. Selon nous, les positions précédentes illustrent deux manières de penser radicalement différentes, qui ont leurs fondements dans les mots eux-mêmes[1]. L'estimation des hommes d'Athènes, l'acte l'emportant sur la parole, n'est pas nécessairement la bonne.

## III. Deux systèmes de pensées.

### 1. Catastrophisme *versus* Évolutionnisme.

Sur l'histoire des concepts, deux interprétations s'affrontent. L'une insiste sur les nouveautés radicales, nées dans un intervalle de temps extrêmement court. L'autre au contraire met l'accent sur son aspect progressif, dont la vitesse dépend certes des périodes considérées. Mais ce qui semble au premier abord une révolution violente, est une apparence, le travail réel, souterrain et invisible, ayant été fait bien auparavant.

Suivant cette analyse (cf. par exemple Pierre Duhem), l'œuvre de Galilée apparaît comme une suite à celle des penseurs médiévaux, ainsi les 'Oxfordiens' ou 'Mertoniens' (du *Merton College*) comme Swineshead ou Heytesbury, mais aussi les 'Parisiens' comme Nicole Oresme autour des $13^{ème}$ et $14^{ème}$ siècles, donc 3 siècles avant Galilée.

Selon l'autre au contraire (cf. Alexandre Koyré), il faut insister sur la rupture d'avec le passé qu'a représentée Galilée, précisément en tant qu'il est à l'origine de la physique moderne ([KOY], p. 172).

On pourrait penser que cette division est confinée aux analyses **sur** les sciences, car la détermination des causes et des auteurs de découvertes ou d'inventions, est toujours délicate. Et pourtant on la retrouve à l'intérieur même des sciences (naturelles). Ainsi est-elle au cœur de la cosmologie, dans la querelle entre 'créationnistes' (s'appuyant sur le '*Big Bang*') et évolutionnistes pour qui l'univers n'a pas de commencement (et peut-être de fin).

C'est en biologie, paléontologie et géologie qu'elle apparaît avec une force nouvelle, au travers d'un débat sans cesse recommencé, situé à l'intersection de ces trois sciences, la 'disparition des dinosaures'. Après divers retournements, il est actuellement un consensus pour l'attribuer à la chute d'une météorite géante (théorie d'Alvarez), signant la victoire de la discontinuité (cf. [BUF], p. 102). Pourtant, rien n'est moins sûr que l'histoire doive s'achever là.

Cette alternative qui apparaît dans presque toutes les sciences, aussi bien que dans leur histoire, recouvre une opposition plus vaste, que nous appellerons les points de vue '**évolutionniste**' et '**catastrophiste**'. Toutefois, les mathématiques posent un problème spécifique.

### 2. Une deuxième alternative.

---

[1] Contrairement à la liberté totale des définitions selon Pascal ([PAS], sect. 1, p. 577).

Il serait en effet possible que les polémiques, dans les sciences de la nature et leur histoire, soient dues aux imprécisions de ce dont elles traitent, concernant généralement des objets ou des événements éloignés, soit de par leur distance, soit temporellement. Car les questions des origines prêtent toujours à controverse, et pas seulement en sciences. Ainsi l'*Enquête* d'Hérodote s'ouvre-t-elle sur 'la responsabilité de la querelle [entre Perses et Phéniciens]' ([HER], I, 1-2).

La connaissance du passé est indirecte et partielle. Mais des mouvements pendulaires des théories standards sont difficilement conciliables avec une immutabilité de la connaissance mathématique, si comme l'affirme Hermann Hankel, contrairement aux autres sciences, en 'mathématiques seulement chaque génération construit une nouvelle histoire sur une structure ancienne'.

Certes les discordes concernant les origines de théorèmes tel celui de Pythagore, les œuvres de Théétète, voire des livres d'Euclide, sont aussi vives que dans les autres sciences. Toutefois, la considération de l'histoire ou la philosophie des mathématiques ne permet d'en rien conclure pour ces dernières. Car ces querelles ont pour objet les mathématiques ou éventuellement les mathématiciens. Elles s'accordent en cela aux histoires des autres disciplines, et ces querelles relèvent alors de la problématique considérée précédemment.

Mais elles ne concernent pas les **objets** des mathématiques (ou des mathématiciens). En tant qu'unique domaine où l'esprit est 'absolument maître', où il est à la fois 'l'outil et ses créations', leurs objets apparaissent intemporels, dépendant exclusivement de l'intellect des mathématiciens qui 'peuvent tisser à leur guise leur propre univers du discours' ([ATK]). Ainsi certains commentateurs peuvent considérer que le calcul infinitésimal ou intégral était déjà (en germes) dans les travaux d'Archimède.

La difficulté se situe dans la relation entre objets réels et objets des mathématiques. En effet, suivant la conception adoptée, on pourra les penser de manières très différentes. Par exemple, soit comme une sorte de formalisation « logique », les résultats se succédant mécaniquement les uns aux autres, soit au contraire à la manière d'une construction architecturale, où la liberté de création est illimitée. Contre Russell, c'est rappelle C. Boyer, l'attitude adoptée par Poincaré pour qui 'si les mathématiciens avaient été la proie de la logique abstraite, ils n'auraient jamais été au-delà de la théorie des nombres et des postulats de la géométrie' ('*Poincaré has said that had mathematicians been left the prey of abstract logic, they would never have gotten beyond the theory of numbers and the postulates of geometry*') ([BOY], p. 13)).

Et tout aussi bien l'existence de différents courants tels les intuitionnistes, les logicistes ou encore les platoniciens[2], qui selon les cas pourraient être classés comme 'évolutionnistes' ou 'catastrophistes', ne permet pas de conclure pour les mathématiques. En effet, ces courants concernent moins les mathématiques que les disciplines qui s'en occupent, histoire ou philosophie. Si une même alternative est possible en mathématiques et dans les autres sciences, elle sera différente de celles que nous avons considérées. Mais ce pourrait en être une généralisation.

### 3. Le discret et le continu.

Si la plus grande partie de cette étude concernant la conception ancienne du mouvement est fondée sur l'œuvre d'Aristote, l'opposition entre discret et continu, reconnue comme problème « philosophique » majeur, peut se lire dans un texte plus « archaïque », si du moins on en croie le Stagirite, la source fondamentale des erreurs de Platon étant sa manière *archaïque* de poser les problèmes ('τό ἀπορῆσαι ἀρχαιχῶς') ([ARI], N, 1089a1).

---

[2] Au sens de la tradition platonicienne en mathématiques.

> '*Socrate (S.) -* **Les choses qui ne sollicitent pas l'intelligence, dis-je, sont celles qui ne suscitent pas simultanément une perception contraire** ; *celles qui suscitent une perception contraire, je considère qu'elles sollicitent l'intelligence, puisque alors leur perception ne manifeste pas plus la chose que ce qui lui est opposé (…) disons que nous avons là trois doigts (…) car jamais la vue ne lui a signifié simultanément qu'un doigt était le contraire d'un doigt.(…) Par conséquent, dis-je, une perception de ce genre ne serait vraisemblablement pas susceptible de solliciter ni d'éveiller l'intelligence.(…) Mais dis-moi, leur grandeur et leur petitesse, la vue les voit-elle de manière satisfaisante ? Et est-ce que cela ne fait aucune différence pour elle que l'un d'entre eux soit placé au milieu ou aux extrémités ? N'en va-t-il pas de même pour le toucher, quand il s'agit de grosseur ou de minceur, de mollesse ou de dureté ? (…) Par conséquent, dis-je, il est nécessaire que dans les cas de ce genre l'âme soit* **perplexe** *(…) et qu'elle se demande aussi, pour la sensation du léger et du lourd, ce qu'est le léger, et le lourd, et* **si la sensation signifie le lourd comme léger, et le léger comme lourd** *?*' (*République*, VII, 523c-524a, traduction un peu modifiée ; c'est nous qui soulignons)[3].

Paradoxalement, Socrate prend ici parti pour le continu contre le discret, pour la quantité contre le nombre ou encore pour la conception géométrique contre la conception arithmétique. En effet, alors que nombrer (en entiers naturels (non nuls)) est un acte purement spontané où l'intelligence n'a point de part, la mesure de ce qui est grand ou petit, gros ou mince, mou ou dur est difficile car engendrant des sensations contradictoires.

Une telle affirmation paraît étrange, puisque la tradition « platonicienne » pose une hiérarchie des sciences où la plus élevée est aussi la plus éloignée de toute application pratique. On a d'ailleurs beaucoup critiqué le philosophe athénien d'avoir interdit le développement des sciences de la nature (physique ou autres), par le mépris dans lequel il tenait toute pratique et en particulier les expériences, parallèlement à Aristote, qui aurait quant à lui, empêché le développement d'une mathématisation du mouvement.

Or la géométrie, ce que rappelle son étymologie, en relation avec la mesure de la terre, est plus proche de nos sens que les nombres qui sont de l'ordre du pur discours (λόγος). Suivant l'école pythagoricienne, que l'on a souvent rapprochée des platoniciens en ce qui concerne les mathématiques, celle-ci tenait que 'tout était nombre', et certainement pas géométrie. Au fondement de toute chose se trouvait la 'τετρακτύς' (la somme des 4 premiers nombres, à savoir 10), 'formidable, toute-puissante, le commencement et le guide de la vie divine aussi bien que terrestre' (Philolaos (-5[ème] siècle e.c.)). Plus généralement on rapporte que Philolaos affirmait que les seules choses connaissables sont des nombres, car 'sans nombre rien ne peut être ni conçu, ni connu', thème que l'on retrouve d'ailleurs beaucoup plus tard, par exemple chez Kronecker, Kummer et même Poincaré.

Il est vrai qu'on ne connaît aucun texte conservé de la période de Pythagore et très peu de son école pendant plusieurs siècles. On doit donc s'en remettre à des commentateurs plus tardifs et engagés, ainsi Aristote. C'est pourquoi il faut comprendre ici Pythagore et pythagorisme comme étant une sorte de tradition.

---

[3] Lorsque Galilée (Salviati) parle de l'esprit sublime de Copernic, ce n'est aucunement pour sa théorie astronomique qui s'accorderait mieux à ce qu'on connaissait du ciel. C'est bien au contraire, parce que contre toutes les données de la perception sensible (αἴσθησις), il a fait prévaloir son intelligence (νόησις) (VII, p. 367). Bien que dans cette étude, il ne soit pas question des démêlés de Galilée avec l'Église, on peut remarquer que la proposition de compromis du cardinal Bellarmin de parler 'par hypothèses' ('*ex suppositione*'), ce qui 'suffit pour le mathématicien' désireux de sauver au mieux 'toutes les apparences' (lettre à Paolo Foscarini, du 12 avril 1615, XII), revient, pour Galilée, à exiger que l'intelligence s'incline devant la sensation (cf. [OFM1], 1[ère] partie).

Pour s'en tenir à la *République*, Platon un peu plus loin place la science des nombres (λογιστική τε καὶ ἀριθμητική) (VII, 525a) avant la géométrie. Mais précise Socrate, cette classification ne suit pas directement de la valeur des objets mathématiques. Elle est fonction inverse du nombre de « dimensions », allant de une pour les 'lignes' (en termes modernes les 'courbes') à trois pour les 'corps' (les 'volumes'), les nombres étant les objets des lignes.

L'alternative catastrophisme/évolutionnisme s'interprète naturellement comme un cas particulier d'opposition du discret et du continu, mais cette fois, les mathématiques participent de cette alternative. Plus encore, elle peut y être posée pour la même raison que dans les autres disciplines. Une connaissance vague et incomplète qui se traduit par des propriétés incohérentes ou inexplicables.

C'est bien ainsi qu'apparaît la quantité irrationnelle, compagne obligée du théorème de Pythagore. Les divers mythes qui marquent sa naissance, son traitement mathématique en rapport à l'infini et son nom lui-même, témoignent des difficultés rencontrées par ceux habitués jusqu'alors à considérer les seuls rapports entiers. C'est aussi l'exemple fondamental du fondateur de l'Académie, sur lequel il ne se lasse pas de revenir.

Un exemple plus récent serait la 'mesure' de Dirac ou la 'fonction' d'Heaviside, 'fonctions' hautement discontinues introduites par les physiciens qui cependant, n'hésitaient pas à les 'dériver'.

La problématique platonicienne permet ainsi de donner un cadre unifié aux mathématiques et aux autres disciplines, en dépit de l'extrême difficulté, du rapport des mathématiques à la réalité. Et, point essentiel, sans qu'il soit nécessaire de se prononcer sur la controverse concernant (la réalité de) ses objets.

**4. Le conflit pour l'éternité.**

Pourtant, comme on l'a dit, que la géométrie soit supérieure à l'arithmétique va à l'encontre de la conception platonicienne usuelle des mathématiques. Et Socrate de poursuivre aussitôt :

'*S. - Dans de tels cas, dis-je, il est vraisemblable que l'âme essaiera en premier lieu, en sollicitant le raisonnement et l'intelligence, d'examiner si chacune des qualités rapportées est unique, ou si elle est double.(...) Par conséquent, s'il apparaît qu'il s'agit de deux choses, c'est que chacune paraît à la fois **différente et une** ? (...) Mais alors, **le nombre et l'unité** [τὸ ἕν], de quel type te semblent-ils ? (...) [Si] **l'unité** peut être vue de manière suffisante (...) elle ne saurait être en mesure de nous tirer vers l'être (...). Si, par contre, on voit toujours simultanément en elle une certaine contradiction, de sorte qu'elle ne semble pas plus être **une** que le contraire (...) l'âme serait nécessairement perplexe et forcée dans ce cas de faire une recherche. Elle mettrait alors en elle-même la réflexion en **mouvement** [κινοῦσα], et elle se demanderait nécessairement ce que peut bien être **l'unité en elle-même.** (...)*
*- Nous voyons en effet [dit Glaucon] simultanément la **même chose comme une et comme une quantité infinie de choses**.*' ([PLA1], VII, 524b-525b, nous soulignons).

Ainsi, de manière très surprenante, Socrate affirme ici le caractère mobile de la réflexion elle-même. Mais en outre, la distinction fondamentale entre discret et continu, que l'on trouvait dans le texte précédent, devient celle entre unité et multiple. On serait sans doute tenté d'associer le discret au multiple, au moins dans le cas fini, mais l'unité participe également du discret. N'est-ce pas contradictoire, si ce n'est avec la première partie du texte platonicien, du moins avec l'analyse que nous en avons donnée au paragraphe précédent ?

Une réponse nécessite une analyse plus approfondie de la notion de 'continu', le recours à des textes modernes mathématiques s'avérant indispensable. Dans *la science et l'hypothèse*, Poincaré rappelle 'la célèbre formule, le continu est l'unité dans la multiplicité' pour ajouter aussitôt que la construction dont se contente l'analyste, intercalant entre des 'échelons' d'autres 'échelons', n'est pas sans poser une difficulté. En effet, le 'continu ainsi conçu n'est plus qu'une collection d'individus (…) *extérieurs* les uns aux autres'. Le continu est certes défini à partir de sa notion physique, identifiée à la conception empirique, mais dont les contradictions ont été supprimées par la puissance symbolique, les mathématiques étant en quelque sorte l'art d'éliminer les contradictions. Quoique mathématiquement correcte, cette définition n'en oublie pas moins l'unité au profit de la multiplicité ([POI1], 1$^{ère}$ partie, chap. II).

Ceci permet de formuler, dans un vocabulaire moderne, la relation entre les deux alternatives du dialogue platonicien. Le continu et donc le discret apparaissent tous deux comme des cas particuliers à la fois de l'unité et de la multiplicité. Ainsi, la topologie moderne rejoint et peut nous aider, nous autres modernes, à comprendre l'argumentation de la *République*.

Si ni le discret ni le continu n'appartiennent exclusivement à l'un ou au multiple, ils en participent des deux ; la contradiction entre les deux parties du dialogue s'estompe donc. Il ne faudrait toutefois pas en déduire un accord entre Platon et Poincaré. Pour le philosophe grec, la raison de distinguer le sens philosophique du continu (et donc du discret) de son sens mathématique, ne peut résulter de ce que l'idée philosophique serait entachée de contradictions.

La seconde alternative enveloppe la première, mais de manière plus complexe qu'une simple généralisation terme à terme. On passe d'un cadre plutôt mathématique (discret et continu) à un cadre beaucoup plus général, où la mathématique est certes toujours présente. C'est la mise en place d'une problématique essentielle, celle de l'un et du multiple.

Celle-ci tient une position centrale dans la pensée platonicienne. On la retrouve dans le *Sophiste* où est énoncée cette autre opposition fondamentale, celle de l'être et du non-être. Il met en scène 'une sorte de combat de géants' (246a), non entre « évolutionnistes » et « catastrophistes », mais 'ceux qui tirent toute chose vers la terre' (246a) (les 'γηγενεῖς' ou 'fils de la terre' (248b)) et les 'amis des idées' ('τῶν εἰδῶν φίλοι' (248a)) 'protégés des hauteurs de l'invisible' (246b). Cette partie du *Sophiste* se confronte à la question du mouvement. Les 'amis des idées', en effet, excluent de la réalité, ce qui est en mouvement (κινεῖσθαι), car seul est, ce qui est au repos (ἠρεμοῦν, de ἠρεμία : tranquillité, repos, terme qu'Aristote utilisera également).

## 5. Du mouvement en mathématiques.

Lorsque Platon résume cette lutte 'sans fin [ἄπλετος : insatiable], [qui] existe à jamais', le problème fondamental n'est pas l'opposition corps/idée. C'est celle du mouvement absolu/repos absolu, chaque partie formant une 'doctrine terrible' (δεινὸν λόγον), puisqu'elles conduisent toutes deux à la 'fin de la science, de la sagesse ou de l'intellect', et à leur propre destruction.

La question du mouvement, plus précisément le point de vue que l'on adopte à son égard, constitue pour le philosophe athénien un problème philosophique majeur, il s'agit de prendre parti dans un combat pour l'éternité. Les positions sont là encore ordonnées entre deux extrêmes, d'une part celle où tout est (en) mouvement, rattachée à la théorie du flux universel héraclitéen, d'autre part celle d'un univers absolument immobile d'où le mouvement, et donc ajoute Platon, la vie elle-même, est exclue.

Il reste que cette lutte de géants entre les partisans du mouvement et ceux du repos n'est pas du même ordre que les alternatives précédentes. Certes, on pourrait bien établir un lien entre la première alternative et le mouvement/repos, encore que l'évolution, même plus lente que la catastrophe, ne soit pas l'immobilité. Pourtant, selon Platon, ce qu'étudient les mathématiques concerne l'immuable et éternel. On pourrait dire que le triangle est le même, qu'il soit considéré par Euclide, nous ou des mathématiciens indiens. D'où un accord général sur les démonstrations mathématiques, mais aussi un accroissement, générations après générations, de connaissances à propos de ce qui est indépendant, à la fois du temps et de l'espace (cf. §III.2). Mais si le mouvement est exclu des mathématiques, aucune alternative incluant celui-là, n'aura de sens pour celles-ci.

Lorsque l'éternité s'ajoute à l'immobilité, on passe du mouvement au changement. Dans la pensée grecque ancienne, les deux concepts étaient conjoints, le mouvement étant, pour reprendre les termes d'Aristote, le changement suivant le lieu. C'est aussi ce que soutient l'Étranger du *Sophiste* de Platon (248a-249e). Nous les modernes, établissons même une équivalence, déjà présente chez Aristote, entre les deux. La 'mobilité', un synonyme de ce qui est apte à se (ou à être) modifier, est généralement conçue positivement, opposant êtres dynamiques à des balourds condamnés à disparaître (à la manière des dinosaures) (M. Caveing souligne la connotation 'axiologique' de la rapidité chez Aristote, cf. [CAV], p. 167-168).

Les objets ou les concepts des mathématiques se situeraient alors du côté des objets fossiles. Pourtant, même si les objets mathématiques ne sont pas considérés comme mobiles ou changeants, le mouvement ou le changement lui-même peut être objet des mathématiques. Il s'agirait d'obtenir des invariants dans ce qui ne l'est pas. Classifier une multitude consiste à trouver certaines caractéristiques, aussi peu nombreuses que possibles, permettant de déterminer les diverses classes. Ainsi en est-il de la classification des courbes analytiques ou des triangles. De même, on pourrait chercher à classifier le mouvement ou le changement au moyen d'invariants qui ne seraient soumis ni à l'un ni à l'autre.

Pour Platon dans le *Sophiste* (mais aussi la *République,* cf. VII, 524b-525b), les mathématiques, aussi bien que les sciences de la nature, n'excluent pas le changement/mouvement. Contrairement à la vulgate platonicienne, il est même impossible de l'éliminer absolument de quelque connaissance que ce soit. Si les 'fils de la terre' sont 'terribles' ('Ἦ δεινοὺς (εἴρηκας) ἄνδρας'), une théorie oublieuse du changement/mouvement le serait également (249a-b). En effet, l'Étranger d'Élée, l'enquêteur du *Sophiste*, affirme qu'il n'est aucune connaissance sans mouvement, non seulement chez celui qui apprend, mais dans la chose qui est en train d'être connue.

On pourrait certes objecter qu'admettre le changement/mouvement en tant qu'objet mathématique, est différent de concevoir, avec l'Étranger du *Sophiste*, que les objets mathématiques eux-mêmes soient sujets à celui-ci. Toutefois, que l'on suive l'hypothèse des invariants ou celle de l'Éléate, il s'ensuit qu'à l'intérieur même des mathématiques, est posée la question du mouvement et du repos.

Il reste qu'on peut se demander si hypothèse il y a. Corrélatifs de toute étude du mouvement, les calculs différentiel et intégral sont si intrinsèquement liés aux mathématiques contemporaines, qu'on aurait beaucoup de peine à trouver un domaine auquel ils ne participent pas. Pourtant la qualité mathématique des infinitésimaux a longtemps fait problème.

Du côté des mathématiciens, Michel Rolle s'emploie à démontrer en 1701 (séances du 12 et 16 mars) devant de l'Académie (royale) des sciences non seulement leur manque de rigueur, mais aussi qu'ils conduisent à des résultats erronés, car contradictoires avec ceux obtenus classiquement d'après les travaux de Fermat et Hudde.

Il considère (en termes modernes) la courbe d'équation :
$$y - b = (x^2 - 2ax + a^2 - b^2)^{2/3} / a^{1/3}$$
et trouve, d'après les règles des infinitésimaux, un unique maximum/minimum en $x=a$. Or suivant les calculs classiques on en obtient trois, les deux autres étant en $x = a+b$ et $x = a-b$ ; d'où la courbe obtenue dont il donne approximativement le dessin suivant :

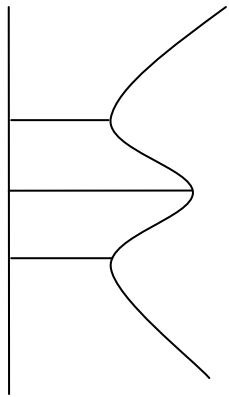

Et Rolle de conclure que l'on 'peut assurer de là que cette géométrie conduit dans l'erreur'.

Dans sa réponse Varignon, partisan des infinitésimaux, montre que Rolle n'a pas appliqué correctement les règles du calcul. Et il corrige la courbe dessinée par Rolle, en montrant que les deux autres points sont certes des minima mais surtout que, contrairement au cas $x = a$, leurs tangentes sont orthogonales à celles dessinées par Rolle. Le lecteur pourra s'amuser à rectifier la figure ci-dessus.

Pour les philosophes, la question reste posée, bien après que les mathématiciens aient intégré le calcul infinitésimal. L'un des plus ingénieux, mais aussi des plus violents adversaires de la théorie des infinitésimaux sous sa forme newtonienne de la théorie des fluxions, est l'évêque George Berkeley qui en 1734 publie '*The Analyst or, a Discourse Addressed to an Infidel Mathematician*' où il entreprend, liant la défense de la religion à la critique des infinitésimaux, de montrer, au mathématicien infidèle (celui qui les défend), les absurdités qu'on y trouve, et ce de la manière la plus impartiale ('*utmost impartiality*', §3). Il faudrait leur concéder, dit-il, des 'entités naissantes et imparfaites' (les flux ou vitesses) mais aussi 'des vitesses de vitesse, les secondes, troisièmes, quatrièmes et cinquièmes vitesses', ce qui dépasse tout entendement humain. Et de s'étonner que des 'libres-penseurs' admettent d'aussi étranges choses ('*most strange*'), eux qui font la fine bouche devant les mystères de la religion : 'Celui qui peut digérer une seconde ou une troisième fluxion ne doit pas être, il me semble, très exigeant au regard de quelque point que ce soit de la divinité' (cf. §7). Mais c'est surtout la propriété 'd'évanescence' qu'il considère comme une absurdité mathématique (§17).

L'importance historique et intellectuelle de ces critiques est soulignée par l'historien des mathématiques, Florian Cajori. C'est pour lui 'l'événement majeur du siècle dans l'histoire des mathématiques britanniques' qui donnera lieu à une très intéressante polémique avec des mathématiciens comme James Jurin ou Benjamin Robins, défenseurs à la fois

Newton et sa théorie des fluxions. Ainsi que le remarque Koyré, les raisons de ces attaques sont les conséquences théologiques qu'impliquent selon ses contemporains, la physique de Newton ([KOY], p. 265). Et inversement, c'est pour répondre à ces attaques que celui-ci ajoute à sa deuxième édition des *Principia* un 'Scholie général' relatif à ses croyances religieuses (*ib.*, p. 268).

Sans doute serait-il malaisé dans les mathématiques contemporaines de retrouver cette controverse, au moins telle qu'elle eut lieu au début du $18^{ème}$ siècle entre Rolle et Varignon. En ce sens, le mouvement/repos est devenu un objet consensuel des mathématiques. Pourtant la querelle toujours vivante entre 'l'intuitionnisme' (plus ou moins modéré) et 'le formalisme' ou le 'logicisme' (plus ou moins modérés) en relèverait encore, si l'on classe la construction avec le changement/mouvement et le symbole avec l'objet/concept immuable (cf. [POI2], livre II, 2, 4 et [POI3], App., II.2).

### 6. Le retour du continu.

Poincaré revient sans relâche sur la question, pour lui fondamentale, du rapport entre les sens mathématique et philosophique du continu. Si le premier a l'avantage de fournir un cadre exempt des contradictions de l'expérience sensible, il le paie d'un prix considérable, une rupture d'avec l'intuition sensible, c'est-à-dire d'avec la nature elle-même ([POI6], 4). L'existence mathématique étant simplement la non-existence de contradictions (id.), ce qui fait problème est l'existence d'une telle réalité (à savoir 'une fonction') qui réponde à ce que l'on (les physiciens) attend des mathématiques. Cette question, reconnaît Poincaré, reflète la querelle entre les philosophies de Leibniz et de Kant. Elle se traduit pour les mathématiciens par celle entre les réalistes cantoriens et les pragmatistes.

Dans notre analyse de la continuité, l'alternative continu/discret présente une différence d'avec les autres questions que nous avons considérées. En effet, l'une des principales difficultés, consistait à montrer que les mathématiques elles-mêmes répondaient à une problématique qui leur était *a priori* extérieure, que ce soit l'alternative évolutionnisme/catastrophisme ou le problème du mouvement/repos. Or cet obstacle n'existe pas pour le continu et le discret, en tant que leurs définitions trouvent leur origines dans les mathématiques. Mais cette simplicité même pourrait cacher un danger : confondre deux définitions différentes sous un même nom. Car il n'est aucunement prouvé que les sens de ces termes en mathématique et en philosophie soient identiques.

Cela pouvait sembler inutile, le texte platonicien par lequel ils étaient introduits portant précisément sur la grandeur, au travers d'une problématique mathématique. Toutefois ils apparaissaient non pas dans un contexte de réflexion philosophique, mais d'une classification de différentes parties des mathématiques (l'arithmétique, la géométrie, l'astronomie, la musique,...). Inversement, l'exemple donné par Socrate ne portait pas directement sur des éléments géométriques ou algébriques, mais sur des perceptions sensibles ; d'une part, visuelles, le dénombrement de doigts, d'autre part, tactiles, le gros et le mince, le dur ou le mou.

C'est précisément contre une telle simplicité que Poincaré met en garde, à savoir confondre la définition mathématique et le sens intuitif d'un même terme. Le continu répète-t-il n'est pas la même chose, en tant que compris par les mathématiciens et en tant que conçu à partir de l'expérience sensible. Certes celle-ci est bien en dernière instance la source de celle-là, il n'empêche que le continu mathématique semble oublier la réunion pour se concentrer sur la division (cf. paragraphe précédent).

Pour préciser cette notion, il faut la référer à ce qu'on nomme actuellement la puissance du continu. Définition purement ensembliste conduisant, dit Poincaré, à des exemples aussi bizarres qu'inutiles que s'amusent à donner certains mathématiciens.

*'Autrefois, quand on inventait une fonction nouvelle, c'était en vue de quelque but pratique ; aujourd'hui, on les invente tout exprès pour mettre en défaut les raisonnements de nos pères, et on n'en tirera jamais que cela.'* ([POI2], livre II, II, 5).

Ainsi en est-il pour lui de l'ensemble de Cantor qui conduit au théorème ou paradoxe de Banach-Tarski selon lequel on peut découper une sphère en un nombre fini (que l'on peut prendre égal à 5) de parties puis en réarranger certaines en sorte d'obtenir deux sphères ou encore une sphère de volume double, et par itération (ou même directement) aussi grande que l'on veut. C'est une autre manière d'exprimer l'existence d'ensembles non Lebesgue mesurables, eux-mêmes conséquences de l'axiome du choix.

Poincaré les critique expressément car, pour lui, dimensions et continu sont intimement liés (*Dernière pensées*, III, 2, p. 28). Si dans de nombreux écrits il se préoccupe de cette question et de son aspect dimensionnel (par ex. [POI4], partie II, 3, 3 ; [POI1], I, 2), c'est dans un texte tardif, qu'il expose le plus clairement sa conception.
Le continu à *n* dimensions est défini comme

*'ensemble de* n *quantités susceptibles de varier indépendamment l'une de l'autre et de prendre toutes les valeurs réelles satisfaisant à certaines inégalités'* ([POI5], III, 2).

Et quoique 'irréprochable', elle n'est pas 'satisfaisante' car liée à des coordonnées. En termes modernes elle n'est pas canonique (ou intrinsèque). Il propose donc une définition par récurrence grâce aux coupures de Dedekind, permettant de le définir simultanément en toutes dimensions.
Un continu de dimension un étant une 'courbe fermée' (*ib.*, p. 28-29), un continu $A$ sera ce qui, lorsqu'on lui ôte un continu $B$ de dimension au moins deux de moins que $A$, reste d'un seul bloc, mais est divisé en deux parties disjointes si $B$ est de dimension un de moins que $A$. Comme on le voit, ces propriétés ne sont pas vérifiées si l'on considère non pas des variétés (lisses) mais des espaces ayant des singularités. L'objectif de Poincaré n'est pas d'élaborer ici une théorie formalisée, mais de mettre en lumière les idées fondamentales concernant la question. Ou plus simplement, son idée du continu était ce qu'on nommerait aujourd'hui une variété différentielle.

Mais ce qui est souligné ici, est l'aspect unitaire du continu, physique ou expérimental, *versus* son absence pour les analystes, autre reproche de Poincaré à leur définition du continu. C'est un retour à la notion historique du continu, le sens moderne ayant été élaboré contre les mathématiciens primitifs ('les pères') qui restaient plus proches de la conception intuitive ou physique.

En effet, Aristote distingue trois notions liées au continu moderne : ce qui est consécutif, contigu et continu (s'y ajoutent encore deux autres que nous ne considérerons pas (cf. *Physique*, V, 3, 226b-227b)). Cette interprétation restera dominante jusqu'aux définitions mathématiques 'rigoureuses' des analystes critiqués par Poincaré.

Consécutif : il n'est aucune séparation du même genre entre les extrémités

Contigu : aucune séparation du même genre et en outre les extrémités sont dans un même lieu.
Continu : contigu dont les extrémités sont identiques.

La *Physique* insiste sur l'aspect de généralité décroissante du consécutif au continu, continu impliquant contigu qui implique consécutif, les réciproques étant fausses. Le contigu serait, du point de vue sensible, un collage si parfait que les extrémités en seraient indiscernables. Pour la mathématique moderne, un exemple serait une réunion de deux intervalles du type [A,B[ et [B,C], tandis qu'une représentation du continu serait une réunion d'intervalles [A,B] et [B,C].

Pourtant le philosophe grec propose également comme exemples de ce qui est continu, 'le clouage, le collage, l'assemblage, la greffe'. Cela peut provenir d'une certaine confusion entre des termes proches dont on n'a pas (encore) de définition (mathématique) rigoureuse, mais il est également possible de les interpréter d'une manière plus cohérente, si l'on considère directement l'argument utilisé dans tous les cas.

En effet, dans le texte aristotélicien, la différentiation des termes se fait en caractérisant les extrémités mises en relation. Le raisonnement sur le continu consiste alors à distinguer celui-ci des autres possibilités, en ce que deux choses n'en forment en réalité qu'une.

Est dit continu, ce qui forme un « un ». C'est l'argument essentiel, car il permet de montrer la contradiction entre indivisible et continu, un élément indivisible n'ayant pas de parties donc d'extrémités. D'où l'impossibilité pour un continu, par exemple le temps, d'être formé de tels indivisibles, les instants. Cette étude se conclut sur une polémique contre l'attribution d'existence réelle au point et à l'unité (ou toute autre notion mathématique pour ce qui concerne Aristote). On ne saurait les identifier car deux points peuvent être joints par une courbe alors qu'il n'est aucun intermédiaire entre deux nombres (entiers). On retrouve la distinction classique chez les anciens Grecs, entre les nombres (entiers) discrets et les grandeurs géométriques continues.

En termes modernes, le continu s'oppose alors au discret, comme l'un au multiple. Toutefois cet « un » serait indéfiniment divisible, alors qu'au contraire ce « multiple » ne le serait, qu'un nombre fini de fois. Ceci devra être néanmoins nuancé, comme on le verra dans [OFM2]. La conception ancienne n'est donc pas transposable directement, puisque on peut avoir topologiquement des ensembles infinis et discrets (par exemple les entiers muni de leur topologie habituelle). Cela importe peu dans le cadre aristotélicien, l'infini n'ayant d'existence que potentielle.

On comprend que 'l'unité' perdue par le continu sous sa forme mathématique, soit précisément celle formée par ce qui est indéfiniment divisible. Dans *les Mathématiques et la Logique,* Poincaré affirme que le fondement intuitif de sa propre définition du continu à partir de coupures, se trouve en effet dans ces divisions successives.

Mais alors ce n'est pas tant le continu auquel nous avons affaire, mais à la **connexité**. Il est vrai que ce n'est qu'avec les travaux sur l'*Analysis Sitûs* (i.e. la topologie), qu'une distinction va être opérée. Une figure géométrique était dite continue, dès lors qu'elle était, d'un seul tenant (l'unité), c'est-à-dire en terme moderne, connexe.

La traduction du terme 'συνεχές' employé par les anciens Grecs et rendu par 'continu' signifie d'abord 'qui se tient' 'qui est ininterrompu' ; pour une figure, sa traduction par 'connexe' serait sans doute plus convenable. En outre, une fonction continue est intuitivement conçue comme envoyant une figure connexe sur une autre, transformant une 'unité multiple' en une autre, par opposition aux fonctions qui laissent des 'trous', ainsi les fonctions en escalier : 'Au début, [l'idée de fonction continue] n'était qu'une image sensible, par exemple, celle d'un trait continu tracé à la craie sur un tableau noir.' (*Valeur de la Science*, I, I, 6, p. 15). Toutefois, si la définition (moderne) de la continuité implique celle-là, l'inverse n'est pas exact.

On en a un exemple avec la fonction réelle : $x \rightarrow sin(1/x)$ pour $x \neq 0$ et s'annulant en *0*. L'image de cette application est certes connexe, et même tout intervalle connexe a une image connexe, en effet l'image d'un voisinage de *0*, aussi petit soit-il, contient (en fait est égal à) tout l'intervalle des réels compris entre *-1* et *1*. Pourtant elle n'est pas continue, la fonction étant très oscillante au voisinage de *0*. C'est précisément la formulation moderne de la continuité qui conduit à distinguer ces deux caractéristiques comme le montre l'observation de Darboux de 1874 ([DAR], p. 109).

La définition proposée par Poincaré se fonde également sur la connexité : une figure formant un bloc qui ne peut être divisée en deux par quelque chose dont la dimension est petite. La propriété d'être continue est globale et se caractérise par l'unité de la figure. C'est encore cette unité qui est absente de la construction des réels par empilement successifs, des entiers aux rationnels aux algébriques et enfin aux réels. Au fur et à mesure, sont comblés les espaces laissés par les empilements précédents, sans jamais obtenir la certitude d'en avoir terminé ([POI1], I, 2). La définition du mathématicien français rejoint la notion aristotélicienne du continu.

Le continu chez Poincaré apparaît posséder une double nature, l'une participant des mathématiques, l'autre de la philosophie (et une troisième peut-être de la physique). Certes, les raisons qu'il en donne, fondées sur la relation entre l'expérience sensible et 'l'esprit humain', ne sont sans doute pas celles, ontologiques, d'Aristote. Et au contraire de Platon, pour qui les mathématiques sont une propédeutique à la philosophie, chez Poincaré, elles en sont l'objet essentiel.

### 7. Les mathématiciens et Platon.

On a souvent dit que les mathématiciens sont spontanément platoniciens, en ce sens qu'ils ne croient pas que, ce sur quoi ils travaillent, soit pure invention imaginaire. Ainsi Alain Connes affirme résolument une '*position platonicienne renouvelée par le théorème de Gödel*' (*La Recherche*, 20, août 2005, p. 77). Cet attachement quasi sentimental n'est peut-être pas entièrement étranger à la place de choix qu'accorde le philosophe athénien aux mathématiques. Une légende ne rapporte-t-elle pas qu'au fronton de l'Académie, lieu qu'il fréquentait, était gravée l'inscription '*Que nul n'entre ici qui ne soit géomètre*' ('ΑΓΕΩΜΕΤΡΗΤΟΣ ΜΗΔΕΙΣ ΕΙΣΙΤΩ') (cf. [SAF])?

Pour Platon, les objets mathématiques ne sont pas des pensées complexes, *inventées* par leur auteur, elles ont une existence propre qui est découverte par l'exploration mathématique du monde des Idées. Il est vrai qu'en mathématique, on ne peut certainement pas faire ce que l'on veut avec ce que l'on a. Ce sur quoi on travaille, possède une consistance, une dureté, une résistance tout à fait semblable à celle d'un objet matériel, et très certainement contraire à une pensée imaginative modifiable à volonté.

'*La perception que nous avons de la réalité mathématique fait que celle-ci se manifeste avec une résistance et une cohérence comparables à celles de la réalité extérieure*' dit A. Connes (*ib.*).

De même, Léon Brunschvicg conclut son analyse des travaux sur la continuité des fonctions en notant, qu'elle

'*s'est imposée aux géomètres et aux analystes presque **en dépit d'eux-mêmes**, en dépit de la tradition séculaire qu'ils avaient tendance à prendre pour intuition immédiate.*' ([BRU], IV, p. 339-340 ; nous soulignons).

Et Jean Dieudonné d'affirmer que '*c'est **contraints** par la **nature profonde** (et souvent cachée jusque-là) des objets et relations classiques, que les mathématiciens, entre 1800 et 1930, ont forgé de nouveaux outils « abstraits »*' ([DIE], intr., p. 10 ; nous soulignons).

Même s'il est vrai que l'on peut trouver des 'mathématiciens contemporains qui se jugent « réalistes »', il n'est pas sûr qu'ils soient si nombreux. Une des difficultés principales à laquelle se heurte le réalisme cantorien critiqué par Poincaré, n'est sans doute pas étrangère à la nécessité, comme dans le cadre des Idées platoniciennes, d'accepter leur délocalisation objective ou ontologique (cf. [CH-CO], p. 35). Il répugne aux mathématiciens, comme à tout un chacun, d'admettre une telle 'délocalisation objective', sous quelque forme que ce soit (cf. A. Connes cité ci-dessus). Sans doute beaucoup seraient plutôt portés vers une forme de 'poincaréisme', en ce sens où les mathématiques se retrouveraient dans la nature humaine, par exemple dans la manière dont l'homme peut concevoir le monde. L'important est que les concepts mathématiques forgés par les mathématiciens soient à la fois leurs inventions propres et en même temps communicables aux autres mathématiciens. On est plus proche d'une certaine forme de kantisme (cf. par ex. [ATK], p. 406-407) que du philosophe athénien. Pour quelques exemples récents, on pourra se référer à deux articles parus dans les *Notices de l'AMS* de décembre 2005, ceux de Martin Gardner et Brian Davies (respectivement p. 1345 et p. 1350)).

Quoique, et c'est en cela qu'ils diffèrent de la plupart des non-mathématiciens, avec un grain plus ou moins étendu de 'réalisme cantorien'. Ainsi, la position du biologiste de [CH-CO] n'est pas tant en désaccord d'avec celle du mathématicien, que témoignant d'une certaine forme d'incrédulité. Ce grain de 'réalisme' serait aussi ce grain de folie qui rendrait nécessaire, selon l'exigence platonicienne, d'être d'abord mathématicien avant de philosopher.

Mais identifier ce 'réalisme' au platonisme, réduit la portée de réalité des Idées platoniciennes et les déforme en projetant sur elles, les objets des mathématiques ou des mathématiciens. Le mathématicien 'réaliste' réclame une délocalisation objective, la possibilité pour les objets mathématiques d'avoir une (certaine) réalité sans toutefois être dans un lieu donné. Pour Platon, on a affaire à une délocalisation ontologique, en ce sens où toute réalité est délocalisée (cf. [PLA2], 130e, 138a-b, 145b-e, 162c). D'un intérêt particulier pour une physique contemporaine oscillant du corpuscule relativiste à l'ondulatoire quantique, est le glissement opéré par Parménide, dans l'ouvrage éponyme. Pour réfuter la délocalisation ontologique avancée par Socrate, il recouvre littéralement d'un voile l'image lumineuse proposée par celui-ci (131b-c). Tous deux s'accordent toutefois pour rejeter une localisation des objets mathématiques dans la pensée (132b), qui aurait la faveur de certains biologistes (cf. [CH-CO]).

La conception de Platon est différente de ce 'réalisme' mathématique, et n'est pas sujet à enorgueillir les mathématiciens classés par le philosophe du côté du chasseur. Ainsi, contrairement à Poincaré qui ne s'y est point trompé dans son opposition au platonisme, les objets considérés par le mathématicien existent indépendamment de lui, et son travail consiste à les traquer. Cette chasse est particulière, puisqu'elle n'a pas pour but de tuer (chasses utilisant les techniques du jet), mais de capturer, son instrument étant le filet qui permet de ramener le gibier vivant, en faisant le moins de dégâts possibles (peut-être ce qu'on appelle le 'prix à payer' en mathématique et qu'en bon économe il faut minimiser). Toutefois, de même que les chasseurs ignorant du savoir de la cuisson et des sauces, ramènent leurs proies au cuisinier qui connaît l'art de les apprêter, les mathématiciens ayant la 'moindre parcelle d'intelligence' savent qu'ils ne peuvent se contenter de garder leurs trouvailles pour eux-mêmes ou leurs collègues, qui ne sauraient qu'en faire. Ils doivent au contraire les porter à qui sait les accommoder convenablement, le (mystérieux) διαλεκτικός, celui qui s'occupe de l'argumentation (cf. *Euthydème*, 290c).

Il serait donc erroné de concevoir les innombrables recours à Platon effectués par les scientifiques du 17ème, et d'abord Galilée, au travers du seul désir de renverser la hiérarchie habituelle entre la scolastique et les mathématiques. Ils recherchent bien plutôt cette pensée de la continuité (cf.), et plus généralement de l'unité, telle qu'elle est présente chez l'auteur des *Lois* (cf. §III.4). Cette unité n'est pas une position moyenne située entre deux systèmes extrêmes, mais au contraire, une pensée en lutte avec elle-même (cf. §III.3). C'est dans cette démarche que se reconnaîtra le savant italien, et plus généralement ceux à l'origine de la nouvelle physique.

Il serait également faux de croire que cette unité consisterait, pour nous autres modernes, en un retour nostalgique à l'idéal du savant total de la Renaissance, à construire des ponts reliant des disciplines. Si l'unité est possible, ce ne peut être que parce qu'il existe une unité de nature dans la pensée. Ou encore, parce que la nature, comme la pensent Galilée et ses successeurs, est continue.

C'est suivant une telle perspective que sera étudiée, dans un autre article ([OFM2]), la naissance d'une nouvelle science, celle du mouvement.

## *TEXTES CITÉS*